\documentclass[a4paper,10pt]{article}
\usepackage[a4paper, total={15cm, 23.7cm}]{geometry}
\usepackage{amssymb,amsthm,amsmath,amsfonts,color,tikz}
\usepackage{comment}
\usepackage[normalem]{ulem}

\newtheorem{theorem}{Theorem}
\newtheorem{lemma}[theorem]{Lemma}
\newtheorem{proposition}[theorem]{Proposition}

\newcommand{\R}{\mathbb{R}}
\newcommand{\N}{\mathbb{N}}

\newcommand{\E}{\mathbb{E}}
\newcommand{\Pa}{\mathbb{P}}

\def\cp{\mathrm{cap}\,}

\newcommand{\Om}{\Omega}

\newcommand{\bp}{\begin{proof}}
\newcommand{\ep}{\end{proof}}

\begin{document}
\title{On some isoperimetric inequalities for the Newtonian capacity}

\author{{M. van den Berg} \\
School of Mathematics, University of Bristol\\
Fry Building, Woodland Road\\
Bristol BS8 1UG\\
United Kingdom\\
\texttt{mamvdb@bristol.ac.uk}}

\date{7 May 2024}\maketitle
\vskip1truecm\indent

\begin{abstract}\noindent
Upper bounds are obtained for the Newtonian capacity of compact sets in $\R^d,\,d\ge 3$ in terms of the perimeter of the $r$-parallel neighbourhood of $K$. For compact, convex sets in $\R^d,\,d\ge 3$ with a $C^2$ boundary the Newtonian capacity is bounded from above by $(d-2)M(K)$, where $M(K)>0$ is the integral of the mean curvature over the boundary of $K$ with equality if $K$ is a ball. For compact, convex sets in $\R^d,\,d\ge 3$ with non-empty interior the Newtonian capacity is bounded from above by $\frac{(d-2)P(K)^2}{d|K|}$ with equality if $K$ is a ball. Here $P(K)$ is the perimeter of $K$ and $|K|$ is its measure. A quantitative refinement of the latter inequality in terms of the Fraenkel asymmetry is also obtained. An upper bound is obtained for expected Newtonian capacity of the Wiener sausage in $\R^d,\,d\ge 5$ with radius $\varepsilon$ and time length $t$.

\end{abstract}

\vskip 1truecm
\noindent\textbf{2020 Mathematics Subject Classification:} 49Q10, 49J45, 49J40, 35J25.\\
\textbf{Keywords}: Newtonian capacity, torsional rigidity, measure, perimeter.

%%%%%%%%%%%%%%%%%%%%%%%%%%%%%%%%%%%%%%%%%%%%%%%%%%%%%%%%%%%%%%%%%%%%%%%%
\section{Introduction \label{sec1}}

In this paper we consider maximisation problems involving Newtonian capacity (or logarithmic capacity if $d=2$).
For a compact set $K\subset \R^d,\,d\ge 3$ we recall a definition of its Newtonian capacity $\cp(K)$ \cite[p.293]{LL}:
\begin{equation}\label{e1}
\cp(K)=\inf\Big\{\int_{\R^d}|D u|^2: u\ge {\bf 1}_K, u\in D^1(\R^d)\cap C^0(\R^d)\Big\},
\end{equation}
where $D^1(\R^d)$ is the collection of functions $f:\R^d\rightarrow \R$ with $f\in L^1_\textup{loc}(\R^d)$, $Df\in L^2(\R^d)$, and which vanish at infinity. Here $f$ vanishes at infinity if for all $\varepsilon>0$, $|\{|f|>\varepsilon\}|<\infty$, where $|A|$ denotes the Lebesgue measure of a measurable set $A\subset \R^d$. The indicator function is denoted by $\bf 1_{\cdot}$.

We introduce the following notation. The boundary of $A$ is denoted by $\partial A$, the perimeter by $P(A)$, the closure by $\overline{A}$, the convex hull by $\textup{co}(A)$, and the
$a$-dimensional Hausdorff measure by $\mathcal{H}^a(A)$ \cite[p.61]{EG}. For a non-empty compact set $K$ we denote for $r>0$ its closed $r$-neighbourhood by
\begin{equation*}%\label{e6}
K_r=\{x\in\R^d: d_K(x)\le r\},
\end{equation*}
where
\begin{equation*}%\label{e7}
d_K(x)=\min\{|x-y|:\,y\in K\},\,x\in \R^d,
\end{equation*}
is the distance to $K$ function. We denote by $C(K)$ the set of critical points of $d_K$. The parallel sets $K_r,\,r>0$ have been studied extensively in the literature. See for example \cite{LS, Fu, RSS, RW}, and the references therein. It is known (\cite[Theorem 4.1]{Fu}) that $C(K)$ is a compact, countable subset of $[0,\infty)$. That theorem also implies that $K_r,\,r\in (0,\infty)\setminus C(K)$ is a Lipschitz manifold, and hence
$$P(K_r) =H^{d-1}(\partial K_r),\,r\in [0,\infty)\setminus C(K).$$
Furthermore in \cite[Theorem 3.3]{RSS} it was shown that
$$\frac{d|K_r|}{d r}=H^{d-1}(\partial K_r),\,r\in [0,\infty)\setminus C(K).$$
Finally in \cite[Lemma 2, Lemma 5]{LS} it was shown that $\frac{d|K_r|}{d r}$ is continuous wherever it exists.

If $K$ is compact then $\R^d\setminus K$ is open and consists of a countable union of open components. Since $K$ is bounded there is precisely one unbounded component of its complement, which is denoted by $U_K$. Let $A_K$ be the union of all bounded components of the complement of $K$. Then $A_K$ is open, and $K=\R^d\setminus(A_K\cup U_K)\subset \R^d\setminus U_K:=\tilde{K}$. It is straightforward to show that $\tilde{K}=K\cup A_K$, and that $\cp(K)=\cp(\tilde{K})$.

If $K$ is compact and $\partial K$ is $C^2$, oriented by an outward unit normal vector field, then we denote the mean curvature map by $H:\partial K\rightarrow \R$, and define its integral by
\begin{equation}\label{e8}
M(K)=\int_{\partial K}Hd\mathcal{H}^{d-1}.
\end{equation}

Our main results are the following.
\begin{theorem}\label{the1} Let $K$ be a non-empty compact set in $\R^d,\,d\ge 3$.
\begin{enumerate}
\item[\textup{(i)}] If  $\int_{(0,\infty)} (P(\tilde{K}_t))^{-1}\,dt<\infty$, then
\begin{equation}\label{e9}
\cp(K)\le \Big(\int_{(0,\infty)} (P(\tilde{K}_t))^{-1}\,dt\Big)^{-1},
\end{equation}
with equality if $K$ is a closed ball.
\item[\textup{(ii)}]If
\begin{equation*}%\label{e11}
\lim_{s\downarrow 0}\int_{(s,\infty)} (P(\tilde{K}_t))^{-1}\,dt=+\infty,
\end{equation*}
then $\cp(K)=0$.
\item[\textup{(iii)}]
\begin{equation}\label{e10}
\cp(K)\le \inf_{a>0}\frac{1}{a^2}|K_a|.
\end{equation}
\item[\textup{(iv)}] If $K$ is convex, and if $\partial K$ is $C^2$, then
\begin{equation}\label{e12}
\cp(K)\le (d-2)M(K),
\end{equation}
with equality if $K$ is a closed ball.
\end{enumerate}
\end{theorem}
It follows from \eqref{e12} and the Aleksandrov-Fenchel inequalities \eqref{e51}, \eqref{e47} below (for $k=2,j=1,i=0$), that
\begin{equation}\label{e13}
\cp(K)\le \frac{(d-2)}{d}\frac{P(K)^2}{|K|},
\end{equation}
with equality if $K$ is any closed ball.

Theorem \ref{the2} below weakens the hypotheses under (iv), and quantifies \eqref{e13} in terms of the Fraenkel asymmetry. The latter is a measure of how close $K$ is to a ball of the same measure as $K$.
For a measurable set $\Om\subset\R^d$ with $0<|\Om|<\infty$ the Fraenkel asymmetry of $\Om$ is the number
\begin{equation}\label{e14}
\mathcal{A}(\Omega)=\inf\Big\{\frac{|\Om\Delta B|}{|B|}: B\, \textup{is a ball with }|B|=|\Om| \Big\}.
\end{equation}
Note that $0\le \mathcal{A}(\Omega)<2$ and that $\mathcal{A}(\Omega)=0$ if and only if  $\Om$ is a ball modulo a set of measure $0$.

It was shown in \cite{BF1} and \cite{BF2} that for $d=2,3,...$ there exist constants $c_d>0$ such that for any compact, convex set $K\subset \R^d$ with $|K|>0$,
\begin{equation}\label{e17}
\frac{P(K)|K|^{-(d-1)/d}}{d\omega_d^{1/d}}-1\ge c_d\mathcal{A}^2(K).
\end{equation}

\smallskip

\begin{theorem}\label{the2} If $K$ is compact and convex in $\R^d,\,d\ge 3$ with $|K|>0$ then
\begin{equation*}%\label{e15}
1-\frac{d\,\cp(K)|K|}{(d-2)P(K)^2}\ge \gamma_d\mathcal{A}^2(K),
\end{equation*}
where
\begin{equation*}%\label{e16}
\gamma_d= \frac{\Gamma(d+1)\Gamma(d-1)}{\Gamma(2d-2)+\Gamma(d)\Gamma(d-1)}\cdot\frac{c_d}{1+4dc_d}\mathcal{A}^2(K).
\end{equation*}
\end{theorem}

In \cite{vdBB} and \cite{MvdBM} the authors obtain inequalities involving  the Newtonian capacity and the torsional rigidity. Recall
that the torsion function for a non-empty open set $\Omega\subset \R^d,\,d\ge 1$ with finite Lebesgue measure is the solution of
\begin{equation}\label{e18}
-\Delta v=1,\quad v\in H_0^1(\Omega),
\end{equation}
and is denoted by $v_{\Omega}$. It is convenient to extend $v_{\Omega}$ to all of $\R^d$ by putting $v_{\Omega}=0$ on $\R^d\setminus\Omega$. The torsion function is non-negative and bounded. Moreover if $\Omega_1,\Omega_2$ are open sets in $\R^d$, then
\begin{equation}\label{e19}
\Omega_1\subset \Omega_2 \Rightarrow 0\le v_{\Omega_1}\le v_{\Omega_2}.
\end{equation}
The torsional rigidity of $\Om$ (or torsion for short) is denoted by
$T(\Omega)=\int_{\Omega}v_{\Omega}$.
Hence \eqref{e19} implies that
\begin{equation}\label{e20}
\Om_1\subset \Om_2 \Rightarrow 0 <T(\Om_1)\le T(\Om_2).
\end{equation}
By \eqref{e18} and the definition of $T(\Om)$,
\begin{equation}\label{e21}
T(t\Om)=t^{d+2}T(\Om),\,t>0,
\end{equation}
where $t\Om$ is a homothety of $\Om$ by a factor $t$.
The de Saint-Venant's inequality \cite[p.206]{AH} asserts that
\begin{equation}\label{e22}
T(\Om)\le T(\Om^*),
\end{equation}
where $\Om^*$ is any ball in $\R^d$ with $|\Om|=|\Om^*|$. By \eqref{e22} and scaling of Lebesgue measure,
\begin{equation}\label{e23}
\frac{T(\Om)}{|\Om|^{(d+2)/d}}\le \frac{T(B_1)}{|B_1|^{(d+2)/d}}=\big(d(d+2)\omega_d^{2/d}\big)^{-1},
\end{equation}
where $B_1$ is the open ball with radius $1$ and measure $\omega_d$.

If $K,K_1,K_2$ are compact sets then
\begin{equation}\label{e2}
K_1\subset K_2 \Rightarrow \cp(K_1)\le \cp(K_2),
\end{equation}
and
\begin{equation}\label{e3}
\cp(tK) =t^{d-2}\cp(K),\,t>0.
\end{equation}
The classical isocapacitary lower bound
\begin{equation*}%\label{e4}
\cp(K)\ge \cp(K^*),
\end{equation*}
where $K^*$ is a closed ball with $|K|=|K^*|$, goes back to \cite{PSZ}. It follows that
\begin{equation}\label{e5}
\frac{\cp(K)}{|K|^{(d-2)/d}}\ge \frac{\cp(\overline{B_1})}{|\overline{B_1}|^{(d-2)/d}}=(d-2)d\omega_d^{2/d},
\end{equation}
where $\overline{B_1}$ is the closed ball with radius $1$.

Let
\begin{equation}\label{e24}
G(\Om)=\frac{T(\Om)\cp(\overline{\Om})}{|\Om|^2}.
\end{equation}
The functional in \eqref{e24} is, by \eqref{e3} and \eqref{e21}, scaling invariant. Fixing $|\Om|=1$ we see that the capacity and torsion are competing: the torsion is, by \eqref{e22}, maximised by a ball with measure $1$, whereas the Newtonian capacity is, by \eqref{e5}, minimised for a ball with measure $1$. Furthermore both torsion and Newtonian capacity are by \eqref{e20} and \eqref{e2} increasing set functions under inclusion.
While the torsion is defined by a variational problem on $\Om$, the Newtonian capacity of $\overline{\Om}$ is defined by a variational problem on its complement. These facts make the study of the variational problems involving $G(\Om)$ very different from the ones leading to the Faber-Krahn inequality or the Kohler-Jobin inequality for example.

In \cite[Theorem 2(i)]{vdBB} it was shown that $G(\Om)$ is not bounded from above on the class of non-empty open sets with finite measure, and in \cite[Theorem 3(i), $q=1$]{vdBB} it was shown that there exists a sequence of convex sets $(\Om_j)$ with $\lim_{j\rightarrow\infty} G(\Om_j)=0$. So the only remaining case of interest is the maximisation of $G(\Om)$ over the collection of convex sets.

In \cite[Theorem 2(iii)]{vdBB} it was shown that the variational problem
\begin{equation}\label{e25}
\sup\{G(\Om) :\, \Om\textup{ non-empty, open, bounded and convex in $\R^d$}\}
\end{equation}
has a maximiser for $d=3$. The existence of a maximiser of the variational problem in \eqref{e25} for $d>3$ remains an open problem.
It was shown in \cite[Theorem 2(i)]{MvdBM} that for any ellipsoid $E\subset\R^d$, $G(E)\le G(B_1)$. This suggests that for any non-empty, open bounded convex set $\Om$, $G(\Om)\le G(B_1)$.

Recall that
\begin{equation}\label{e27}
P(t\Om)=t^{d-1}P(\Om),\, t>0.
\end{equation}
For $d\ge 3$ and $0\le \alpha\le 2$ we define the functional
\begin{equation}\label{e28}
G_{\alpha}(\Om)=\frac{T(\Om)\cp(\overline{\Om})}{|\Om|^{\alpha}P(\Om)^{d(2-\alpha)/(d-1)}}.
\end{equation}
By \eqref{e3}, \eqref{e21} and \eqref{e27}, we see that $G_{\alpha}$ is scaling invariant. The functional interpolates between a perimeter and a measure constraint.
The following was shown in \cite[Theorem 6]{MvdBM}:
\begin{enumerate}
\item[\textup{(i)}]Let $\mathfrak E_d$ denote the collection of open ellipsoids in $\R^d$. If $d\ge 3$ and $0\le \alpha\le2$, then
\begin{equation}\label{e29}
\sup\{G_{\alpha}(\Omega): \Omega\in \mathfrak{E}_d\}=G_{\alpha}(B_1),
\end{equation}
and the supremum in the left-hand side of \eqref{e29} is achieved if and only if $\Omega$ is a ball.
\item[\textup{(ii)}] If $d\ge 3$ and $0\le \alpha\le2$, then
\begin{equation}\label{e30}
 \sup\{G_{\alpha}(\Om):\Om\,\, \textup{non-empty, open, bounded, convex in $\R^d$}\}  \le d^{2d}G_{\alpha}(B_1).
\end{equation}
\item[\textup{(iii)}]
If $0\le\alpha<2$, then the  variational problem in the left-hand side of \eqref{e30} has a maximiser.
\end{enumerate}
The presence of a perimeter term in the denominator of  $G_{\alpha}$ guarantees the existence of a maximiser.

 \smallskip

Theorem \ref{the3} below shows that $B_1$ is a maximiser for $G_{\alpha}$ among the collection of open bounded convex sets provided the exponent of the perimeter is not too small. Theorem \ref{the3} together with (i) and (iii) above suggest that $B_1$ is a maximiser of \eqref{e28} for $0\le\alpha\le2$.

\begin{theorem}\label{the3}
If $d\ge 3$ and $0\le \alpha\le \frac{2}{d}$, then
\begin{equation*}%\label{e32}
\sup\{G_{\alpha}(\Omega): \Omega\,\, \textup{non-empty, open, bounded, convex in}\, \R^d\}=G_{\alpha}(B_1),
\end{equation*}
and any ball is a maximiser of $G_{\alpha}$.
\end{theorem}

This paper is organised as follows. In Section \ref{sec2} below we prove Theorems \ref{the1}, \ref{the2} and \ref{the3}. Section \ref{sec3} concerns the analysis of some variational problems involving collections of open sets $\Om\subset \R^2$ with torsion $T(\Om)$, and logarithmic capacity $\cp(\overline{\Om})$. In Section \ref{sec4} we give various examples and discuss the optimality for the bounds in Theorems \ref{the1} and \ref{the3}.

\section{Proofs of Theorems \ref{the1}, \ref{the2} and \ref{the3} \label{sec2}}

\medskip

\noindent {\it Proof of Theorem \ref{the1}.} The starting point of the proof of Theorem \ref{the1} goes back to Theorem 11 in \cite{CFG} where the authors obtain, for convex bodies $K$, an upper bound for $\cp(K)$ by restricting the test functions in \eqref{e1} to those depending on $d_K$ only.
The proof of Theorem \ref{the1}(i) is organised as follows. In step (a) we restrict the class of test functions in \eqref{e1} and derive, formally, a candidate for a test function. In step (b) we show that this function is well defined.
In steps (c)--(e) we show that this function satisfies the constraints in \eqref{e1}, and is admissible. We then complete the proof of (i).

\smallskip

\noindent(a) Let $s>0$ be arbitrary, and let $\varphi=f(d_{\tilde{K}_s})$, where $\tilde{K}=\R^d\setminus U_K$. Then
\begin{equation*}%\label{e33}
|D\varphi|^2=(f'(d_{\tilde{K}_s}))^2|Dd_{\tilde{K}_s}|^2\le (f'(d_{\tilde{K}_s}))^2.
\end{equation*}
By the coarea formula
\begin{equation}\label{e34}
\int_{\R^d}|D\varphi|^2\le \int_{\R^d}(f'(d_{\tilde{K}_s}))^2=\int_{(0,\infty)} (f'(r))^2P(\tilde{K}_{s+r})\,dr.
\end{equation}
Minimising, formally, over all smooth $f$ with $f(0)=1, f(\infty)=0$ gives $(f'(r)P(\tilde{K}_{s+r}))'=0.$ Hence $f'(r)P(\tilde{K}_{s+r})=c$ for some $c\in \R$. It follows that
\begin{equation}\label{e35}
f(r)=c\int_{(r,\infty)}(P(\tilde{K}_{s+t}))^{-1}\,dt=c\int_{(r+s,\infty)}(P(\tilde{K}_{t}))^{-1}\,dt.
\end{equation}
Since $f(0)=1$ we find that
\begin{equation}\label{e38}
f(r)=\frac{\int_{(r,\infty)} (P(\tilde{K}_{s+t}))^{-1}\,dt}{\int_{(0,\infty)}(P(\tilde{K}_{s+t}))^{-1}\,dt}.
\end{equation}

\smallskip

\noindent(b) Since $K\ne \emptyset$ it contains a point say $0$. Then $\overline{B(0;t)}\subset\tilde{K}_t$. Since $\overline{B(0;t)}$ is a convex subset of $\tilde{K}_{t}$,
\begin{equation}\label{e36}
P(\tilde{K}_{t})\ge P(B_1)t^{d-1},
\end{equation}
and
\begin{equation}\label{e36a}
P(\tilde{K}_{t})\ge P(B_1)(r+s)^{d-1},\,t\ge r+s.
\end{equation}
Since $C(\tilde{K})$ is compact, $(r+s,\infty)\setminus C(\tilde{K})$ is open, and hence is a countable union of disjoint open intervals. By the properties of parallel sets mentioned in Section \ref{sec1}, $t\mapsto P(\tilde{K}_{t})$ is continuous on $(r+s,\infty)\setminus C(\tilde{K})$. By \eqref{e36a}, $P(\tilde{K}_{t})$ is uniformly bounded away from $0$. Hence  $t\mapsto (P(\tilde{K}_{t}))^{-1}$ is continuous on $(r+s,\infty)\setminus C(\tilde{K}).$
Since $C(\tilde{K})$ is countable it has measure $0$.
This shows that the integral in the right-hand side of \eqref{e35} is well-defined.
To show that the integral in \eqref{e35} converges we have, by \eqref{e36}, that
\begin{equation*}%\label{e37}
\int_{(r,\infty)}(P(\tilde{K}_{s+t}))^{-1}dt\,\le ((d-2)P(B_1)(s+r)^{d-2})^{-1}\le((d-2)P(B_1)s^{d-2})^{-1}<\infty.
\end{equation*}

\smallskip

\noindent(c) To prove continuity we have by \eqref{e38} and \eqref{e36}.
\begin{align*}%\label{e39}
|\varphi(x)-\varphi(y)|&\le |f(d_{\tilde{K}_s}(x))-f(d_{\tilde{K}_s}(y))|\nonumber\\& \le \sup_{r\ge 0}|f'(r)||d_{\tilde{K}_s}(x)-d_{\tilde{K}_s}(y)|\nonumber\\&\le
\frac{\sup_{r\ge 0}(P(\tilde{K}_{s+r}))^{-1}}{\int_{(0,\infty)} (P(\tilde{K}_{s+t}))^{-1}\,dt}|x-y|\nonumber\\&
\le\frac{(d\omega_ds^{d-1})^{-1}}{\int_{(0,\infty)} (P(\tilde{K}_{s+t}))^{-1}\,dt}|x-y|.
\end{align*}
Hence $\varphi$ is uniformly continuous. This in turn implies that $f\in L^1_\textup{loc}(\R^d)$.

\smallskip

\noindent(d) To prove that $\varphi$ vanishes at infinity in the sense of $D^1(\R^d)$, we have by \eqref{e36} and \eqref{e38}
\begin{equation}\label{e40}
f(r)\le \frac{\big((d-2)d\omega_dr^{d-2}\big)^{-1}}{\int_{(0,\infty)} (P(\tilde{K}_{s+t}))^{-1}\,dt}.
\end{equation}
By \eqref{e40}
\begin{align}\label{e41}
\{\varphi>\varepsilon\}&=\{x\in\R^d: f(d_{\tilde{K}_s}(x))>\varepsilon\}\nonumber\\&\subset\Big\{x\in\R^d: d_{\tilde{K}_s}(x)<\Big(\varepsilon(d-2)d\omega_d \int_{(0,\infty)} (P(\tilde{K}_{s+t}))^{-1}dt\,\Big)^{-1/(d-2)}\Big\}.
\end{align}
Since $\tilde{K}_s$ is contained in a ball with radius $\textup{diam}(\tilde{K}_s),$ we have by \eqref{e41} that the level set $\{\varphi>\varepsilon\}$
is contained in a ball with radius
\begin{equation*}%\label{e42}
\textup{diam}(\tilde{K}_s)+\Big((d-2)d\omega_d\varepsilon \int_{(0,\infty)} (P(\tilde{K}_{s+t}))^{-1}\,dt\Big)^{-1/(d-2)}.
\end{equation*}
Hence this level set has finite measure.

\smallskip

\noindent(e) To see that $D\varphi\in L^2(\R^d)$ we compute
by \eqref{e34} and \eqref{e38} that
\begin{equation}\label{e43}
\int_{\R^d}|D\varphi|^2\le\Big(\int_{(0,\infty)} (P(\tilde{K}_{s+t}))^{-1}\,dt\Big)^{-1}<\infty.
\end{equation}
We conclude by (c)--(e) above that $\varphi\in D^1(\R^d)\cap C(\R^d)$, and hence is a test function. By \eqref{e1} and \eqref{e43},
\begin{equation}\label{e44}
\cp(K)=\cp(\tilde{K})\le \cp(\tilde{K}_s)\le \Big(\int_{(0,\infty)} (P(\tilde{K}_{s+t}))^{-1}\,dt\Big)^{-1}<\infty.
\end{equation}
Since $s>0$ was arbitrary we arrive at \eqref{e9}. A direct computation gives equality for a ball. This completes the proof of Theorem \ref{the1} (i).

\smallskip

\noindent (ii) The assertion follows immediately from \eqref{e44}.

\smallskip

\noindent (iii) Let $a>0$ be arbitrary. Since $P(\tilde{K}_t)\le P(K_t),$ we have by \eqref{e9},
\begin{equation}\label{e9a}
\cp(K)\le \Big(\int_{(0,\infty)} (P(K_t))^{-1}\,dt\Big)^{-1}.
\end{equation}
By Cauchy-Schwarz and \eqref{e9a}
\begin{align*}%\label{e45}
a=\int_{(0,a)} dt&\le \Big(\int_{(0,a)}(P(K_t))^{-1}dt\,\Big)^{1/2}\Big(\int_{(0,a)}P(K_t)\,dt\,\Big)^{1/2}\nonumber\\&\le
\Big(\int_{(0,\infty)}(P(K_t))^{-1}\,dt\Big)^{1/2}|K_a|^{1/2}\nonumber\\&
\le \cp(K)^{-1/2}|K_a|^{1/2}.
\end{align*}
This implies \eqref{e10} since $a>0$ was arbitrary.

\smallskip

\noindent (iv) Steiner's formula for non-empty, compact and convex $K$ reads
\begin{equation}\label{e46}
|K_r|=\sum_{n=0}^d \binom{d}{n} W_n(K)r^n,
\end{equation}
where the $W_n(K)$ are the Quermass integrals for $K$. See \cite[Sections (4.1), (4.2) Chapter 4]{RS}.
These Quermass integrals can be expressed in terms of integrals over the surface $\partial K$ of polynomials in the $d-1$ principal curvatures.
In particular
\begin{equation}\label{e47}
W_0(K)=|K|,\,\, W_1(K)=d^{-1}P(K),\,\, W_2(K)=d^{-1}\int_{\partial K} Hd\mathcal{H}^{d-1},\,\, W_d(K)=\omega_d.
\end{equation}
The right-hand side of \eqref{e46} is differentiable. Hence
\begin{equation}\label{e48}
P(K_r)=\frac{d|K_r|}{d r}=\sum_{n=1}^d n\binom{d}{n} W_n(K)r^{n-1}.
\end{equation}
By the change of variable
\begin{equation*}%\label{e49}
r=\frac{W_1(K)}{W_2(K)}\theta,
\end{equation*}
we obtain by \eqref{e9} and \eqref{e48}
\begin{equation}\label{e50}
\int_{(0,\infty)} (P(K_r))^{-1}dr\,\ge \frac{1}{W_2(K)}\int_{(0,\infty)} \Big(\sum_{n=1}^d n\binom{d}{n} \frac{W_n(K)W_1(K)^{n-2}}{W_2(K)^{n-1}}\theta^{n-1}\Big)^{-1}d\theta.
\end{equation}
The Aleksandrov-Fenchel inequalities \cite[(7.66)]{RS} read
\begin{equation}\label{e51}
W_j(K)^{k-i}\ge W_i(K)^{k-j}W_k(K)^{j-i},\,\,0\le i<j<k\le d.
\end{equation}
Let $j=2,\,i=1,\,k=n$ in \eqref{e51}. This gives
\begin{equation}\label{e52}
W_n(K)W_1(K)^{n-2}\le W_2(K)^{n-1}.
\end{equation}
By \eqref{e50} and \eqref{e52},
\begin{align*}%\label{e53}
\int_{(0,\infty)} (P(K_r))^{-1}dr\,&\ge\frac{1}{W_2(K)}\int_{(0,\infty)} \Big(\sum_{n=1}^d n\binom{d}{n} \theta^{n-1}\Big)^{-1}d\theta\nonumber\\&
=\frac{1}{d(d-2)W_2(K)}\nonumber\\&
=\frac{1}{(d-2)M(K)},
\end{align*}
where we have used \eqref{e8} and \eqref{e47}. A direct computation gives equality for a ball. This proves (iv) by \eqref{e9}, and completes the proof of Theorem \ref{the1}.
\hspace*{\fill }$\square $

\medskip

\noindent {\it Proof of Theorem \ref{the2}.} Let $r>0$. By Steiner's formula \eqref{e46} applied to the compact, convex set $K_r$,
\begin{equation*}%\label{e54}
|K_{r+s}|=\sum_{n=0}^d \binom{d}{n} W_n(K_r)s^n,\,s>0,
\end{equation*}
where the Quermass integrals $W_n(K_r)$ satisfy the Aleksandrov-Fenchel inequalities. By the change of variable
\begin{equation*}%\label{e55}
s=\frac{W_0(K_r)}{W_1(K_r)}\theta,
\end{equation*}
we obtain
\begin{equation}\label{e56}
\int_{(0,\infty)}\frac{ds}{P(K_{r+s})}=\frac{W_0(K_r)}{W_1(K_r)^2}\int_{(0,\infty)}\bigg(\sum_{n=1}^d n\binom{d}{n} \frac{W_n(K_r)W_0(K_r)^{n-1}}{W_1(K_r)^n}\theta^{n-1}\bigg)^{-1}d\theta.
\end{equation}
We put $j=1,i=0,k=n$ in \eqref{e51} to get that
\begin{equation*}%\label{e57}
W_1(K)^{n}\ge W_0(K)^{n-1}W_n(K),\,\,1\le n\le d.
\end{equation*}
This together with \eqref{e9}, \eqref{e47} and  \eqref{e56} gives
\begin{align}\label{e58}
\frac{P(K_r)^2}{\cp(K_r)|K_r|}&\ge d^2\int_{(0,\infty)} \bigg(\sum_{n=1}^{d-1} n\binom{d}{n} \theta^{n-1}+d\frac{W_d(K_r)W_0(K_r)^{d-1}}{W_1(K_r)^d}\theta^{d-1}\bigg)^{-1}d\theta\nonumber\\&
=d\int_{(0,\infty)}\bigg((1+\theta)^{d-1}-\theta^{d-1}\bigg(1-\frac{W_d(K_r)W_0(K_r)^{d-1}}{W_1(K_r)^d}\bigg)\bigg)^{-1}d\theta\,.
\end{align}
Since the integrand in the first line of \eqref{e58} is positive, we have that
\begin{equation}\label{e59}
\theta^{d-1}\bigg(1-\frac{W_d(K_r)W_0(K_r)^{d-1}}{W_1(K_r)^d}\bigg)<(1+\theta)^{d-1}.
\end{equation}
By \eqref{e59}
\begin{align}\label{e60}
\frac{P(K_r)^2}{\cp(K_r)|K_r|}&\ge d\int_{(0,\infty)}(1+\theta)^{1-d}\bigg(1-\frac{\theta^{d-1}}{(1+\theta)^{d-1}}\bigg(1-\frac{W_d(K_r)W_0(K_r)^{d-1}}{W_1(K_r)^d}\bigg)\bigg)^{-1}d\theta\,\nonumber\\&\ge
d\int_{(0,\infty)}(1+\theta)^{1-d}\bigg(1+\frac{\theta^{d-1}}{(1+\theta)^{d-1}}\bigg(1-\frac{W_d(K_r)W_0(K_r)^{d-1}}{W_1(K_r)^d}\bigg)\bigg)d\theta\,\nonumber\\&=
\frac{d}{d-2}+\frac{\Gamma(d+1)\Gamma(d-2)}{\Gamma(2d-2)}\bigg(1-\frac{W_d(K_r)W_0(K_r)^{d-1}}{W_1(K_r)^d}\bigg),
\end{align}
where we have used \cite[3.194.3]{GR}.

By \eqref{e47} and \eqref{e17}
\begin{align}\label{e61}
\frac{W_d(K_r)W_0(K_r)^{d-1}}{W_1(K_r)^d}&=\bigg(d\omega_d^{1/d}|K_r|^{(d-1)/d}P(K_r)^{-1}\bigg)^d\nonumber\\&
\le\bigg(\frac{1}{1+c_d\mathcal{A}^2(K_r)}\bigg)^d\nonumber\\&\le
\frac{1}{1+dc_d\mathcal{A}^2(K_r)}\nonumber\\&\le 1-\frac{dc_d\mathcal{A}^2(K_r)}{1+4dc_d},
\end{align}
where we have used that the Fraenkel asymmetry is bounded from above by $2$. By \eqref{e60} and \eqref{e61},
\begin{equation}\label{e62}
\frac{P(K_r)^2}{\cp(K_r)|K_r|}\ge
\frac{d}{d-2}+\frac{\Gamma(d+1)\Gamma(d-2)}{\Gamma(2d-2)}\frac{dc_d}{1+4dc_d}\mathcal{A}^2(K_r).
\end{equation}
Rewriting \eqref{e62} as
\begin{equation}\label{e63}
1-\frac{d\,\cp(K_r)|K_r|}{(d-2)P(K_r)^2}\ge\frac{C}{1+C}
\end{equation}
gives
\begin{align}\label{e64}
C&=\frac{\Gamma(d)\Gamma(d-1)}{\Gamma(2d-2)}\frac{dc_d}{1+4dc_d}\mathcal{A}^2(K_r)\nonumber\\ &\le
\frac{\Gamma(d)\Gamma(d-1)}{\Gamma(2d-2)}.
\end{align}
By \eqref{e63} and \eqref{e64}
\begin{equation}\label{e65}
1-\frac{d\,\cp(K_r)|K_r|}{(d-2)P(K_r)^2}\ge \frac{\Gamma(d+1)\Gamma(d-1)}{\Gamma(2d-2)+\Gamma(d)\Gamma(d-1)}\cdot\frac{c_d}{1+4dc_d}\mathcal{A}^2(K_r).
\end{equation}
By monotonicity $\cp(K_r)\ge \cp(K)$, and $|K_r|\ge |K|$. So \eqref{e65} implies
\begin{equation*}%\label{e65a}
1-\frac{d\,\cp(K)|K|}{(d-2)P(K_r)^2}\ge \frac{\Gamma(d+1)\Gamma(d-1)}{\Gamma(2d-2)+\Gamma(d)\Gamma(d-1)}\cdot\frac{c_d}{1+4dc_d}\mathcal{A}^2(K_r).
\end{equation*}
By \cite[2.4.1--2.4.3]{BB} we have that $\lim_{r\downarrow 0}P(K_r)=P(K).$ It therefore suffices to show the following.
\begin{lemma}\label{lem1} If $K$ is a compact, convex set with non-empty interior, then
\begin{equation}\label{e65b}
\lim_{r\downarrow 0}\mathcal{A}(K_r)=\mathcal{A}(K).
\end{equation}
\end{lemma}
The straightforward proof is included for completeness.
\begin{proof} If $\mathcal{A}(K)=0$ then $K$ is a ball, and so is $K_r$. Then $\mathcal{A}(K_r)=0$, and there is nothing to prove.
Suppose $\mathcal{A}(K)>0$. To prove the lemma we let $B_r$ be the ball which minimises the right-hand side of \eqref{e14} with $\Omega=K_r$,
and let $0$ be its centre. We denote by $B_0$ the ball with that same centre $0$, and measure $|K|$. We have by \eqref{e14},
\begin{align}\label{e65i}
\mathcal{A}(K)&\le \frac{|K\Delta B_0|}{|B_0|}
\le\frac{|K_r\Delta B_0|}{|B_0|}+\frac{|K_r\setminus K|}{|B_0|}\nonumber\\&
\le\frac{|K_r\Delta B_r|}{|B_0|}+\frac{|K_r\setminus K|}{|B_0|}+\frac{|B_r\setminus B_0|}{|B_0|}\nonumber\\&
=\mathcal{A}(K_r)+\mathcal{A}(K_r)\Big(\frac{|B_r|}{|B_0|}-1\Big)+\frac{|K_r\setminus K|}{|B_0|}+\frac{|B_r\setminus B_0|}{|B_0|}\nonumber\\&
\le \mathcal{A}(K_r)+\frac{4|K_r\setminus K|}{|K|},
\end{align}
where we have used that
$|K_r|=|B_r|,|K|=|B_0|$, and that $\mathcal{A}(K_r)\le 2$ in the last line of \eqref{e65i}.

We now let $B_0'$ be the ball which minimises the right-hand side of \eqref{e14} with $\Omega=K$, and let $0$ be its centre.
We denote by $B_r'$ the ball with that same centre $0$ and measure $|K_r|$. We have by \eqref{e14}
\begin{align}\label{e65j}
\mathcal{A}(K_r)&\le \frac{|K_r\Delta B_r'|}{|B_r'|}
\le \frac{|K_r\Delta B_r'|}{|B_0'|}\nonumber\\&
\le \frac{|K\Delta B_r'|}{|B_0'|}+\frac{|K_r\setminus K|}{|K|}\nonumber\\&
\le \frac{|K\Delta B_0'|}{|B_0'|}+\frac{|K_r\setminus K|}{|K|}+\frac{|B_r'\setminus B_0'|}{|B_0'|}\nonumber\\&
=\mathcal{A}(K)+\frac{2|K_r\setminus K|}{|K|}.
\end{align}
By \eqref{e65i} and \eqref{e65j} we find that
\begin{equation}\label{e65k}
|\mathcal{A}(K_r)-\mathcal{A}(K)|\le \frac{4|K_r\setminus K|}{|K|}.
\end{equation}
By \cite[2.4.1--2.4.3]{BB} we have that $\lim_{r\downarrow 0}|K_r|=|K|$. This, together with \eqref{e65k}, gives \eqref{e65b}.
\end{proof}
This completes the proof of Theorem \ref{the2}.
\hspace*{\fill }$\square $

\medskip

\noindent {\it Proof of Theorem \ref{the3}.}
By Theorem \ref{the2} and definition \eqref{e28}
\begin{align}\label{e66}
G_{\alpha}(\Om)&\le \frac{d-2}{d}\frac{T(\Om)}{|\Om|^{(d+2)/d}}\frac{|\Om|^{(2-\alpha d)/d}}{P(\Om)^{(2-\alpha d)/(d-1)}}\nonumber\\ &
\le \frac{d-2}{d}\frac{T(B_1)}{|B_1|^{(d+2)/d}}\frac{|B_1|^{(2-\alpha d)/d}}{P(B_1)^{(2-\alpha d)/(d-1)}}\nonumber\\ &
=G_{\alpha}(B_1).
\end{align}
We have used the de Saint-Venant's inequality \eqref{e23} for the first fraction in the right-hand side of \eqref{e66}, the isoperimetric inequality for the second fraction, and definition \eqref{e28} for the last equality.
\hspace*{\fill }$\square $

\section{Logarithmic capacity \label{sec3}}

In this section we denote by $\cp(\cdot)$ the logarithmic capacity, defined on the class of compact sets in $\R^2$, and recall its definition below.
Let $\mu$ be a probability measure supported on $K$, and let
\begin{equation*}%\label{e67}
I(\mu)=\iint_{K\times K}\log\Big(\frac{1}{|x-y|}\Big)\mu(dx)\mu(dy).
\end{equation*}
Furthermore let
\begin{equation*}%\label{e68}
V(K)=\inf\big\{I(\mu):\mu \textup{ a probability measure on $K$} \big\}.
\end{equation*}
The logarithmic capacity of $K$ is denoted by $\cp(K)$, and is the non-negative real number
$\cp(K)=e^{-V(K)}.$

The logarithmic capacity is an increasing set function, and satisfies \eqref{e2} for compact sets $K_1$ and $K_2$ in $\R^2$.
For an ellipsoid with semi-axes $a_1$ and $a_2$,
\begin{equation*}%\label{e69}
\cp(\overline{E(a)})=\frac12(a_1+a_2).
\end{equation*}
See \cite{L}.

Let $d=2,\,0\le \alpha\le \frac32 $, and let
\begin{equation*}%\label{e70}
H_{\alpha}(\Omega)=\frac{T(\Omega)^{1/2}\cp(\overline{\Omega})}{|\Omega|^{\alpha}P(\Om)^{3-2\alpha}}.
\end{equation*}
Then $H_{\alpha}$ is scaling invariant.
The following results were obtained in \cite[Theorem 7]{MvdBM}:
\begin{enumerate}
\item[\textup{(i)}]Let $\mathfrak E_2$ denote the collection of open ellipses in $\R^2$. If $0\le \alpha\le \frac{3}{2}$, then
\begin{equation}\label{e71}
\sup\{H_{\alpha}(\Omega): \Omega\in \mathfrak{E}_2\}=H_{\alpha}(B_1),
\end{equation}
and the supremum in the left-hand side of \eqref{e71} is achieved if and only if $\Omega$ is a ball.
\item[\textup{(ii)}]If $0\le\alpha\le\frac32$, then
\begin{equation}\label{e72}
 \sup\{H_{\alpha}(\Om):\Om\, \textup{non-empty, open, bounded, convex}\}\le 2^{2\alpha}\pi^{3-2\alpha}H_{\alpha}(B_1).
\end{equation}
\item[\textup{(iii)}] If $0\le\alpha<\frac32$, then the variational problem in the left-hand side of \eqref{e72} has a maximiser.
If $\Omega_{\alpha}$ is any such maximiser, then
\begin{equation}\label{e73}
\frac{\textup{diam}(\Omega_{\alpha})}{\rho(\Omega_{\alpha})}\le 2^{(3+2\alpha)/(3-2\alpha)}\pi^2,
\end{equation}
where $\rho(\cdot)$ denotes the inradius.
\item[\textup{(iv)}]If $\alpha=0$, then the variational problem
\begin{equation*}%\label{e74}
\sup\big\{H_{0}(\Omega):\Omega\,\textup{\,open, bounded, connected, $0<|\Om|<\infty$}\big \},
\end{equation*}
has a maximiser. Any such maximiser is also a maximiser of \eqref{e72} for $\alpha=0$, and henceforth satisfies \eqref{e73}.
\end{enumerate}

The main result of this section is the following.
\begin{theorem}\label{the4} If $d=2$, then
\begin{equation}\label{e75}
\sup\Big\{\frac{T(\Om)\cp(\overline{\Om})}{P(\Om)^{5}}: \Om\,\textup{non-empty, open, bounded, connected}\Big\}=\frac{T(B_1)\cp(\overline{B_1})}{P(B_1)^{5}},
\end{equation}
and $B_1$ is a maximiser of the left-hand side of \eqref{e75}.
\end{theorem}
\begin{proof} We have
\begin{align}\label{e76}
\frac{T(B_1)\cp(\overline{B_1})}{P(B_1)^{5}}&\le\sup\Big\{\frac{T(\Om)\cp(\overline{\Om})}{P(\Om)^{5}}: \Om\,\textup{non-empty, open, bounded, connected}\Big\}\nonumber\\&
=\sup\Big\{\frac{\cp(\overline{\Om})}{P(\Om)}\frac{|\Om|^2}{P(\Om)^4}\frac{T(\Om)}{|\Om|^2}: \Om\,\textup{non-empty, open, bounded, connected}\Big\}\nonumber\\&
\le \sup\Big\{\frac{\cp(\overline{\Om})}{P(\Om)}\frac{|B_1|^2}{P(B_1)^4}\frac{T(B_1)}{|B_1|^2}: \Om\,\textup{non-empty, open, bounded, connected}\Big\},
\end{align}
where we have used in the final inequality in \eqref{e76} the isoperimetric inequality, and the de Saint-Venant's inequality respectively.
It is clear that $\Omega$ is contained in the closure of its convex hull $\overline{\textup{co}(\Om)}$. Hence  $\cp(\overline{\Om})\le \cp(\overline{\textup{co}(\Om)})$.
Furthermore since $\Omega$ is connected $P(\Om)\ge P(\overline{\textup{co}(\Om)})$. By inequality \cite[Table 1.21, Formula 12]{PSZ} we have for any bounded convex set $A\subset \R^2$,
\begin{equation}\label{e77}
\frac{\cp(\overline{A})}{P(A)}\le \frac{\cp(\overline{B_1})}{P(B_1)}.
\end{equation}
Applying \eqref{e77} to the convex set $\overline{\textup{co}(\Om)}$, and using \eqref{e76} we arrive at \eqref{e75}.
\end{proof}

\section{Examples and Optimality \label{sec4}}

The example below shows that there exist compact sets $K\subset\R^3$ with $\cp(K)=0$ for which the right-hand side of \eqref{e9} is strictly positive. It is straightforward to find such examples for $d>3$.
\noindent \begin{proposition}\label{prop1}
Let $\alpha>0$, let $n\in \N$, and let $K(\alpha)\subset \R^3$ be given by
\begin{equation*}%\label{e78}
K(\alpha)=\Big(\bigcup_{n\in\N}\{(n^{-\alpha},0)\}\cup\{(0,0)\}\Big)\times[0,1].
\end{equation*}
\begin{itemize}
\item[\textup{(i)}] $K(\alpha)$ is a compact subset of $[0,1]^3$ with $\cp(K(\alpha))=0$.
\item[\textup{(ii)}] If $\alpha>0$ then
\begin{equation}\label{e79}
\Big(\int_{(0,\infty)}(P(K(\alpha)_r))^{-1}dr\,\Big)^{-1}\ge \frac{4\pi\alpha}{2^{\alpha+2}+3\alpha(\alpha+1)}.
\end{equation}
\end{itemize}
\end{proposition}
\begin{proof}
\noindent(i) Since $K(\alpha)$ is a countable union of line segments in $\R^3$, $\cp(K(\alpha))=0.$

\smallskip

\noindent(ii) Let
\begin{equation}\label{e80}
r^*=\frac{\alpha}{2^{\alpha+2}}.
\end{equation}
We wish to obtain a lower bound for $P(K(\alpha)_r).$ For $r\ge r^*$ we use that $P(K(\alpha)_r)\ge 4\pi r^2$,
and find by \eqref{e80} that
 \begin{equation}\label{e81}
\int_{(r^*,\infty)}(P(K(\alpha)_r))^{-1}dr\,\le \frac{1}{4\pi r^*}=\frac{2^{\alpha}}{\pi\alpha}.
\end{equation}
To obtain a lower bound for $P(K(\alpha)_r)$  for $0<r\le r^*$ we consider all pairs of line segments which are at least distance $2r$ apart.
The distance between line segments with  $x_1=n^{-\alpha}$ and $x_1=(n+1)^{-\alpha}$ is bounded from below by
 \begin{equation*}%\label{e82}
n^{-\alpha}-(n+1)^{-\alpha}\ge \alpha(n+1)^{-\alpha-1}.
\end{equation*}
So if $2r\le \alpha(n+1)^{-\alpha-1}$ then all line segments with $x_1\ge n^{-\alpha}$ contribute at least $2\pi r$ to the perimeter.
There are at least $n_r$ of such line segments, where
\begin{equation*}%\label{e83}
n_r= \Big[\Big(\frac{\alpha}{2r}\Big)^{1/(\alpha+1)}\Big]-1,
\end{equation*}
and where $[\cdot]$ denotes the integer part. Hence
\begin{align*}%\label{e84}
P(K(\alpha)_r)&\ge 2\pi r\Big(\Big[\Big(\frac{\alpha}{2r}\Big)^{1/(\alpha+1)}\Big]-1\Big).
\end{align*}
For all $x\ge 2$ we have $[x]-1\ge \frac{x}{3}.$ The choice in \eqref{e80} implies that
\begin{equation*}%\label{e85}
\Big[\Big(\frac{\alpha}{2r}\Big)^{1/(\alpha+1)}\Big]-1\ge \frac13 \Big(\frac{\alpha}{2r}\Big)^{1/(\alpha+1)},\,0\le r\le r^*.
\end{equation*}
Hence
\begin{equation*}%\label{e86}
P(K(\alpha)_r)\ge \frac{2\pi}{3}\Big(\frac{\alpha}{2}\Big)^{1/(\alpha+1)}r^{\frac{\alpha}{\alpha+1}},\,0<r\le r^*,
\end{equation*}
and by \eqref{e80}
\begin{align}\label{e87}
\int_{(0,r^*)}(P(K(\alpha)_r))^{-1}dr\,&\le\frac{3}{2\pi}\Big(\frac{2}{\alpha}\Big)^{1/(\alpha+1)}(1+\alpha)(r^*)^{1/(1+\alpha)}\nonumber\\ &=
\frac{3}{4\pi}(1+\alpha).
\end{align}
By \eqref{e81} and \eqref{e87}
\begin{equation*}%\label{e88}
\int_{(0,\infty)}(P(K(\alpha)_r))^{-1}dr\,\le \frac{2^{\alpha}}{\pi\alpha}+\frac{3}{4\pi}(1+\alpha).
\end{equation*}
This implies \eqref{e79}.
\end{proof}

\smallskip

Below we show that the maximisation of $|K|^{\alpha}\cp(K)$ over all compact, convex sets in $\R^d,\,d\ge 3$ with given perimeter leads either to a restatement of \eqref{e13} for $\alpha\ge 1$, or
an infinite supremum for $0<\alpha<1$. So the exponent $1$ of $|K|$ in the variational problem
\begin{equation*}%\label{e89}
\sup\Big\{\frac{|K|\cp(K)}{P(K)^2}: K\,\, \textup{non-empty, compact, convex in}\, \R^d\Big\}
\end{equation*}
is optimal. The statement under \eqref{e13} asserts that $\overline{B_1}$ is a maximiser.

Define for $\alpha>0$ the scaling invariant functional
\begin{equation}\label{e90}
J_{\alpha}(K)=\frac{|K|^{\alpha}\cp({K})}{P(K)^{(d\alpha+d-2)/(d-1)}}.
\end{equation}

\smallskip

\begin{proposition}\label{prop2}
\begin{itemize}
\item[\textup{(i)}] If $d\ge 3$ and $\alpha\ge 1$, then
\begin{equation*}%\label{e91}
\sup\{J_{\alpha}(K): K\,\, \textup{non-empty, compact, convex in}\, \R^d\}= J_{\alpha}(B_1),
\end{equation*}
so that $\overline{B_1}$ is a maximiser of the left-hand side of \eqref{e90}.
\item[\textup{(ii)}] If $d\ge 3$ and $0<\alpha< 1$, then
\begin{equation*}%\label{e92}
\sup\{J_{\alpha}(K): K\,\, \textup{non-empty, compact, convex in}\, \R^d\}=+\infty.
\end{equation*}
\end{itemize}
\end{proposition}
\begin{proof}
To prove (i) we rewrite $J_{\alpha}$ as follows:
\begin{equation}\label{e93}
J_{\alpha}(K)=\frac{|K|\cp(K)}{P(K)^2}\Big(\frac{|K|}{P(K)^{d/(d-1)}}\Big)^{\alpha-1}.
\end{equation}
The first term in the right-hand side of \eqref{e93} is, by Theorem \ref{the2}, bounded for compact, convex sets in $\R^d$ by $\frac{|\overline{B_1}|\cp(\overline{B_1})}{P(\overline{B_1})^2}$.
The second term in the right-hand side of \eqref{e93} is bounded from above by the isoperimetric inequality,  $\Big(\frac{|\overline{B_1}|}{P(\overline{B_1})^{d/(d-1)}}\Big)^{\alpha-1}$. This proves the assertion under (i).

To prove (ii) we consider the open ellipsoid $E_{\varepsilon}$ with  $d-2$ semi-axes of length $1$ and $2$ semi-axes of length $\varepsilon$, where $0<\varepsilon<1$ is arbitrary. We have that
\begin{equation}\label{e94}
|E_{\varepsilon}|=\omega_d\varepsilon^2.
\end{equation}
Since $E_{\varepsilon}$ is contained in the cuboid $(-\varepsilon,\varepsilon)\times(-\varepsilon,\varepsilon)\times(-1,1)^{d-2}$ we have that
\begin{equation}\label{e95}
P(E_{\varepsilon})\le P\big((-\varepsilon,\varepsilon)\times(-\varepsilon,\varepsilon)\times(-1,1)^{d-2}\big)\le d2^d\varepsilon.
\end{equation}

Let $E(a)$, with $a=(a_1,a_2,\dots,a_d)\in\R_+^d$, be the ellipsoid
\begin{equation*}%\label{e96}
E(a)=\bigg\{x\in\R^d\ :\ \sum_{i=1}^d\frac{x_i^2}{a_i^2}<1\bigg\}.
\end{equation*}
It was reported in \cite[p.260]{IMcK} that the Newtonian capacity of an ellipsoid was computed in \cite[Volume 8, p.30]{GC}. The formula there is for a three-dimensional ellipsoid, and is given in terms of an elliptic integral. It extends to all $d\ge 3$, and reads
\begin{equation}\label{e97}
\cp\big(\overline{E(a)}\big)=2d\omega_d\mathfrak{e}(a)^{-1},
\end{equation}
where
\begin{equation}\label{e98}
\mathfrak{e}(a)=\int_0^{\infty}dt\,\bigg(\prod_{i=1}^d\big(a_i^2+t\big)\bigg)^{-1/2}.
\end{equation}
In \cite{MvdBM}, \eqref{e97} and \eqref{e98} were used to obtain upper bounds on the capacity of an ellipsoid. Below we bound $\mathfrak{e}(a)$ from above to obtain a lower bound for $\cp(\overline{E_{\varepsilon}})$.
\begin{align}\label{e99}
\cp(\overline{E_{\varepsilon}})&=2d\omega_d\Big(\int_0^{\infty}(1+t)^{-(d-2)/2}(\varepsilon^2+t)^{-1}dt\,\Big)^{-1}\nonumber\\ &\ge
2d\omega_d\Big(\int_0^{\infty}(1+t)^{-1/2}(\varepsilon^2+t)^{-1}dt\,\Big)^{-1}\nonumber\\ &=
2d\omega_d(1-\varepsilon^2)^{1/2}\Big(\log\Big(\frac{1+(1-\varepsilon^2)^{1/2}}{1-(1-\varepsilon^2)^{1/2}}\Big)\Big)^{-1}\nonumber\\ &\ge
d\omega_d(1-\varepsilon^2)^{1/2}\Big(\log\Big(\frac{2}{\varepsilon}\Big)\Big)^{-1},
\end{align}
where we have used that $d\ge 3$ in the second line of \eqref{e99}. By \eqref{e90}, \eqref{e94}, \eqref{e95} and \eqref{e99} we conclude
\begin{equation*}%\label{e100}
J_{\alpha}(E_{\varepsilon})\ge \frac{d\omega_d^{1+\alpha}}{(d2^d)^{(d\alpha+d-2)/(d-1)}}\varepsilon^{(d-2)(\alpha-1)/(d-1)}\frac{(1-\varepsilon^2)^{1/2}}{\log\big(\frac{2}{\varepsilon}\big)},\,0<\varepsilon<1.
\end{equation*}
Hence $J_{\alpha}(E_{\varepsilon})$ is not bounded from above since $0<\alpha<1$, and $\varepsilon\in (0,1)$ was arbitrary. This proves the assertion under (ii).
\end{proof}

In Proposition \ref{prop3} we obtain some elementary information on the Newtonian capacity of the Wiener sausage for a compact set $K$ in $\R^d$. The notation and construction is as follows.
Let $(\beta(s),s\ge 0; \Pa_x,x\in\R^d)$ be Brownian motion, that is the Markov process with generator $\Delta$. Here $\Pa_x$ is the law of $\beta(\cdot)$ starting at $x$ with corresponding expectation $\E_x$.
The Wiener sausage of (time) length $t$ associated to the compact set $K$ is the random set (\cite{JFG,vdBBdH})
\begin{equation*}%\label{e101}
W_t^{K}=\bigcup_{0\le s\le t}\big(\beta(s)+K\big).
\end{equation*}
Since the Brownian path is continuous a.s. we have that the Wiener sausage up to $t$ is a compact set a.s.

\begin{proposition}\label{prop3}
If $d\ge 5$, and if $K$ is a compact set, then
\begin{enumerate}
\item[\textup{(i)}]
\begin{equation}\label{e102}
\limsup_{t\rightarrow\infty}\frac{1}{t}\E_0(\cp(W_t^{K}))\le 16\inf_{c>0}\frac{1}{c^4}|K_c|.
\end{equation}
\item[\textup{(ii)}] If $d\ge 5$, and if $K=\overline{B_{\varepsilon}}=\varepsilon\overline{B_1},\,\varepsilon>0$, then
\begin{equation*}%\label{e103}
\limsup_{t\rightarrow\infty}\frac{1}{t}\E_0(\cp(W_t^{\overline{B_{\varepsilon}}}))\le\kappa_d\frac{(d-2)^{d-2}}{4(d-4)^{d-4}}\varepsilon^{d-4}.
\end{equation*}
\end{enumerate}
\end{proposition}
In fact the $\limsup_{t\rightarrow\infty}$ in the left-hand side of \eqref{e102} could be replaced by $\lim_{t\rightarrow\infty}$. See \cite[(1.8)]{vdBBdH}.
\begin{proof}
To prove the inequality we use classical results going back to \cite{FS} and to \cite[Theorems 1,\,2,\,3]{JFG}. These imply that for $d>2$,
\begin{equation}\label{e104}
\lim_{t\rightarrow\infty}\frac{1}{t}\E_0(|W_t^{K}|)=\cp(K).
\end{equation}
Since
\begin{equation*}%\label{e105}
(W_t^{K})_a=W_t^{K_a},\, a>0,
\end{equation*}
we have
\begin{equation}\label{e106}
|(W_t^{K})_a|=|W_t^{K_a}|,\, a>0,
\end{equation}
and by \eqref{e10}, \eqref{e104} and \eqref{e106},
\begin{align}\label{e107}
\E_0(\cp(W_t^{K}))&\le \E_0\Big(\frac{1}{a^2}|W_t^{K_a}|\Big)\nonumber\\&
=\frac{1}{a^2}\E_0(|W_t^{K_a}|)\nonumber\\&
=\frac{1}{a^2}\cp(K_a)t(1+o(1)),\,t\rightarrow\infty.
\end{align}
Using \eqref{e10} once more, for the compact set $K_a$, we obtain,
\begin{equation*}%\label{e108}
\E_0(\cp(W_t^{K}))\le \frac{1}{a^2}\frac{1}{b^2}|K_{a+b}|t(1+o(1)),\,t\rightarrow\infty.
\end{equation*}
Choosing $a=b=\frac{c}{2}$ yields the assertion under (i).

To prove (ii) we use that
\begin{equation}\label{e109}
\cp((\overline{B_{\varepsilon}})_a)=\kappa_d(a+\varepsilon)^{d-2}.
\end{equation}
This gives by \eqref{e107} and \eqref{e109},
\begin{equation}\label{e110}
\E_0(\cp(W_t^{B_{\varepsilon}}))\le \frac{1}{a^2}\kappa_d(a+\varepsilon)^{d-2}t(1+o(1)).
\end{equation}
Minimising the right-hand side of \eqref{e110} with respect to $a$ gives for $d\ge 5$ with $a=\frac{2\varepsilon}{d-4}$,
\begin{equation*}%\label{e111}
\E_0(\cp(W_t^{\overline{B_{\varepsilon}}}))\le \kappa_d\frac{(d-2)^{d-2}}{4(d-4)^{d-4}}\varepsilon^{d-4}t(1+o(1)),\, t\rightarrow\infty.
\end{equation*}
This implies the assertion under (ii).
\end{proof}

\textbf{Acknowledgements:} MvdB would like to thank Dorin Bucur, Mar\'ia Hern\'andez Cifre and Jan Rataj for very helpful discussions. He is also grateful for hospitality at the Laboratoire de Math\'ematiques, Universit\'e Savoie Mont Blanc, Le-Bourget-Du-Lac in November 2022.

\medskip

%\noindent\textbf{Funding:}

\end{document}